\documentstyle[12pt,graphicx,pstricks,braket,amsmath,amssymb,enumitem,
                   float, epstopdf,pdfpages,cite,tabu,hyperref,subfigure, float,tikz]{article}
                    \usetikzlibrary{arrows,matrix,positioning}
\begin{document}
	\def\NPB#1#2#3{{\it Nucl.\ Phys.}\/ {\bf B#1} (#2) #3}

\def\AEF{A.E. Faraggi}

\begin{titlepage}
\samepage{
\setcounter{page}{1}
\rightline{}
\vspace{1.5cm}

\begin{center}
 {\Large \bf The Hydrogen Atom, Pi and Lerch's Transcendent}
\end{center}

\begin{center}

{\large
Johar M. Ashfaque$^\spadesuit$
}\\
\vspace{1cm}
$^\spadesuit${\it  Max Planck Institute for Software Systems,
            Campus E1 5,\\ 66123 Saarbr\"ucken, Germany\\}
\end{center}

\begin{abstract}
In this note, we extend the connection between the hydrogen atom and $\pi$ to the number $e$ via the Lerch's transcendent.
\end{abstract}
\smallskip}
\end{titlepage}

\section{Introduction}
Wallis' product for $\pi$ written down in 1655 by John Wallis states that
$$\prod_{n=1}^{\infty}\bigg(\frac{2n}{2n-1}\cdot \frac{2n}{2n+1}\bigg)=\frac{\pi}{2}.$$

\section{Proof Via Euler Infinite Product For The Sine Function}
The Euler infinite product for the sine function is
$$\frac{\sin x}{x} =\prod_{n=1}^{\infty}\bigg(1-\frac{x^{2}}{n^{2}\pi^{2}}\bigg).$$
Letting $x=\frac{\pi}{2}$ which gives
$$\frac{2}{\pi} = \prod_{n=1}^{\infty}\bigg(1-\frac{1}{4n^{2}}\bigg).$$
Then
$$\frac{\pi}{2} = \prod_{n=1}^{\infty} \bigg(\frac{4n^{2}}{4n^{2}-1}\bigg) = \prod_{n=1}^{\infty} \bigg(\frac{(2n) (2n)}{(2n-1)(2n+1)}\bigg).  $$
\section{Wallis' Product \& The Hydrogen Atom}
In \cite{Friedmann:2015, Friedmann:2017vde, Chashchina:2017wzn}, the link between the hydrogen atom and $\pi$ has been explored through various trial functions like the Gaussian trial function \cite{Friedmann:2015, Friedmann:2017vde}, the Lorentz trial function \cite{Chashchina:2017wzn} and finally \cite{ Cortese:2017din} who took one step further and showed that the Wallis formula is related to the harmonic oscillator by using the duality between the hydrogen atom and the harmonic oscillator as well as showing that the asymptotic formula under consideration is related to the ``clever" choice of the trial function and a potential in the Schr\"odinger equation. 

Here, for simplicity, we sketch some of the details in the case of the Gaussian trial function, that can be found in \cite{Friedmann:2015, Friedmann:2017vde, Chashchina:2017wzn, Cortese:2017din}. The main idea is to consider the limit
$$\lim_{l\rightarrow \infty}\frac{\langle H \rangle^{l}_{min}}{E_{0,l}} = \lim_{l\rightarrow \infty} \frac{(l+1)^{2}}{(l+\frac{3}{2})}\bigg[\frac{\Gamma(l+1)}{\Gamma(l+\frac{3}{2})} \bigg]^{2}=1.$$
and then let $n=l+1$ giving
$$\lim_{n\rightarrow \infty} \frac{n^{2}}{(n+\frac{1}{2})}\bigg[\frac{\Gamma(n)}{\Gamma(n+\frac{1}{2})} \bigg]^{2}=1.$$
Thereafter making use of the fact that 
$$\bigg(n+\frac{1}{2}\bigg)\Gamma \bigg(n+\frac{1}{2}\bigg)= \Gamma\bigg(n+\frac{3}{2}\bigg)$$
we obtain
$$\lim_{n\rightarrow \infty} \frac{\Gamma(n+1)^{2}}{\Gamma(n+\frac{1}{2})\Gamma(n+\frac{3}{2})} =1.$$
Now making use of the Legendre duplication formula
$$\Gamma(2n) = \frac{2^{2n-1}\Gamma(n)\Gamma(n+\frac{1}{2})}{\sqrt{\pi}}\Rightarrow \frac{1}{\Gamma(n+\frac{1}{2})} = \frac{2^{2n-1}\Gamma(n)}{\sqrt{\pi}\Gamma(2n)}$$
we find that 
$$\frac{\Gamma (n+1)}{\Gamma(n+\frac{1}{2})} = \frac{1}{\sqrt{\pi}}\prod_{j=1}^{n} \frac{2j}{2j-1} $$
and
$$\frac{\Gamma (n+1)}{\Gamma(n+\frac{3}{2})} = \frac{2}{\sqrt{\pi}}\prod_{j=1}^{n} \frac{2j}{2j+1}$$
Multiplying the last two results gives
$$\frac{\Gamma(n+1)^{2}}{\Gamma(n+\frac{1}{2})\Gamma(n+\frac{3}{2})} = \frac{2}{\pi} \prod_{j=1}^{n} \frac{(2j)(2j)}{(2j-1)(2j+1)} $$
and taking the limit, we arrive at the Wallis' product
$$\lim_{n\rightarrow \infty} \prod_{j=1}^{n} \frac{(2j)(2j)}{(2j-1)(2j+1)} =\frac{\pi}{2}.$$

\section{The Lerch's transcendent and $\pi$}
Following \cite{GuiSondow}, the Lerch's transcendent $\Phi$ is defined to be the analytic continuation of the series
$$\Phi(z, s, u) = \frac{1}{u^s} + \frac{z}{(u+1)^s}+\frac{z^2}{(u+2)^s}+...$$
which converges for any real number $u >0$ if $z$ and $s$ are any complex numbers with either $|z|<1$, or $|z|= 1$ and $\Re(s)>1$. The following identity of the Lerch's transcendent will be of key importance 
$$\Phi(z, s-1, u) = \bigg(u+z\frac{\partial}{\partial z}\bigg)\Phi(z, s, u).$$
The corollary in \cite{GuiSondow} then states that for $m= 0,1,2,...$ and complex $z$ with $\Re(z)<1/2$,
$$\frac{\partial \Phi}{\partial s}(z,-m,u) =\bigg(u+z\frac{\partial}{\partial z} \bigg)^m\sum_{n=0}^{\infty}\frac{1}{1-z}\bigg(\frac{-z}{1-z}\bigg)^n\sum_{n=0}^{n} (-1)^{k+1}\binom{n}{k}\ln(u+k).$$
Then setting $m= 0$, $z=−1$ and $u= 1$ and multiplying by $2$ gives
$$2\frac{\partial \Phi}{\partial s}(-1,0,1) =\sum_{n=0}^{\infty}\frac{1}{2^n}\sum_{k=0}^{n}(-1)^{k+1}\binom{n}{k}\ln(k+1).$$
Now making use of the facts that 
$$B_1(x) = -\Phi(1,0,x) = x-\frac{1}{2}$$
where $B_1$ is the first Bernoulli polynomial and
$$\frac{\partial \Phi}{\partial s}(1,0,u)=\ln \frac{\Gamma(u)}{\sqrt{2\pi}}$$
we have that 
$$\frac{\partial \Phi}{\partial s}(-1,0,u)=\ln \frac{\Gamma(\frac{u}{2})}{\Gamma(\frac{u+1}{2})\sqrt{2}}.$$
This identity yields the desired result that
$$\frac{\pi}{2} =\sqrt{2} \bigg(\frac{2^2}{3}\bigg)^{\frac{1}{4}} \bigg(\frac{2^3\cdot 4}{3^3}\bigg)^{\frac{1}{8}} \bigg(\frac{2^4\cdot 4^4}{3^6\cdot 5}\bigg)^{\frac{1}{16}}...=e^s$$
where 
$$s= 2\frac{\partial \Phi}{\partial s}(-1,0,1).$$

\end{document}